\documentclass{amsart}

\usepackage{mathtools}
\usepackage{mathabx}
\usepackage{booktabs}
\usepackage[T1]{fontenc}
\usepackage{amsmath,amsfonts,mathrsfs,amssymb,amsthm}
\usepackage{aliascnt}
\usepackage{enumerate}
\usepackage{bm}
\usepackage{dsfont}
\usepackage{stmaryrd}
\usepackage{wrapfig}
\usepackage{tikz}
\usepackage{tikz-3dplot}
\definecolor{myblue}{rgb}{0.21, 0.34, 0.74}
\definecolor{mygrey}{rgb}{0.55, 0.57, 0.67}
\definecolor{myred}{rgb}{0.79, 0.0, 0.09}
\definecolor{mygreen}{rgb}{0.05, 0.5, 0.06}

\usepackage[colorlinks=true, citecolor=myblue, linkcolor=black]{hyperref}
\DeclareMathAlphabet{\mathscrbf}{OMS}{mdugm}{b}{n}
\usepackage[scr=boondoxo]{mathalpha}
\usepackage{cleveref}
\usepackage{subcaption}
\DeclareCaptionLabelFormat{empty}{}
\usepackage{comment}

\usepackage{wrapfig}
\usepackage{upgreek}

\renewcommand{\S}{\textsection}

\newcommand{\Id}{\mathrm{Id}}

\newcommand{\proj}{\mathbf{P}}

\DeclareMathOperator{\dive}{\mathring{\mathrm{d{\imath}v}}}

\renewcommand{\S}{{\mathbb{S}^{d-1}}}

\newcommand{\cM}{\mathcal{M}}

\newcommand{\cP}{\mathcal{P}}

\newcommand{\cS}{\mathcal{S}}

\newcommand{\cU}{\mathcal{U}}

\newcommand{\cX}{\mathcal{X}}

\newcommand*\diff{\mathop{}\!\mathrm{d}}
\newcommand{\R}{{\rm I}\kern-0.18em{\rm R}}
\newcommand{\h}{{\rm I}\kern-0.18em{\rm H}}
\newcommand{\K}{{\rm I}\kern-0.18em{\rm K}}
\newcommand{\p}{{\rm I}\kern-0.18em{\rm P}}
\newcommand{\E}{{\rm I}\kern-0.18em{\rm E}}

\newcommand{\1}{{\rm 1}\kern-0.24em{\rm I}}
\newcommand{\N}{{\rm I}\kern-0.18em{\rm N}}

\numberwithin{equation}{section}

\usepackage[font=small,labelfont=bf]{caption}

\theoremstyle{plain}
\newtheorem{theorem}{Theorem}[section]

\newaliascnt{proposition}{theorem}
\newtheorem{proposition}[proposition]{Proposition}
\aliascntresetthe{proposition}
\crefname{proposition}{Proposition}{Propositions}
\Crefname{proposition}{Proposition}{Propositions}

\newaliascnt{lemma}{theorem}
\newtheorem{lemma}[lemma]{Lemma}
\aliascntresetthe{lemma}
\crefname{lemma}{Lemma}{Lemmas}
\Crefname{lemma}{Lemma}{Lemmas}

\newaliascnt{corollary}{theorem}
\newtheorem{corollary}[corollary]{Corollary}
\aliascntresetthe{corollary}
\crefname{corollary}{Corollary}{Corollaries}
\Crefname{corollary}{Corollary}{Corollaries}

\newaliascnt{definition}{theorem}

\aliascntresetthe{definition}
\crefname{definition}{Definition}{Definitions}
\Crefname{definition}{Definition}{Definitions}

\newaliascnt{remark}{theorem}

\aliascntresetthe{remark}
\crefname{remark}{Remark}{Remarks}
\Crefname{remark}{Remark}{Remarks}

\newaliascnt{example}{theorem}

\aliascntresetthe{example}
\crefname{example}{Example}{Examples}
\Crefname{example}{Example}{Examples}

\newaliascnt{problem}{theorem}

\aliascntresetthe{problem}
\crefname{problem}{Problem}{Problems}
\Crefname{problem}{Problem}{Problems}

\newaliascnt{assumption}{theorem}
\newtheorem{assumption}[assumption]{Assumption}
\aliascntresetthe{assumption}
\crefname{assumption}{Assumption}{Assumptions}
\Crefname{assumption}{Assumption}{Assumptions}

\crefname{claim}{Claim}{Claims}
\Crefname{claim}{Claim}{Claims}

\crefname{experiment}{Experiment}{Experiments}
\Crefname{experiment}{Experiment}{Experiments}

\begin{document}

\title[Uniform-in-time propagation of chaos on the sphere]{Uniform-in-time propagation of chaos for the Kuramoto model on the sphere}

\author[S. Chen]{Shi Chen}
\address{(SC) Department of Mathematics, Massachusetts Institute of Technology, 77 Massachusetts Ave, 02139 Cambridge MA, USA} 
\email{schen636@mit.edu}

\author[Z. Lin]{Zhengjiang Lin}
\address{(ZL) Department of Mathematics, Massachusetts Institute of Technology, 77 Massachusetts Ave, 02139 Cambridge MA, USA} 
\email{linzj@mit.edu}

\author[E. Parmaksiz]{Emre Parmaksiz}
\address{(EP) Department of Mathematics, Massachusetts Institute of Technology, 77 Massachusetts Ave, 02139 Cambridge MA, USA} 
\email{yep@mit.edu}

\date{}


\begin{abstract}
We establish quantitative uniform-in-time propagation of chaos for a generalization of the classical Kuramoto model on spheres of arbitrary dimension. 
For $n$ particles that are initialized independently from a probability measure with an $L^p$ density, $p>1$, and a nonzero mean vector, we obtain an explicit high-probability bound for the supremum over all times of the Wasserstein distance between the empirical and mean-field distributions.
In sufficiently high dimension, this bound has the same order as the standard empirical sampling rate on the sphere and is therefore optimal. 
The proof combines a finite-time propagation of chaos estimate, a quantitative synchronization result for the mean-field model, and a late-time stability estimate of the system. 
\end{abstract}
\maketitle

\setcounter{tocdepth}{2}
\makeatletter
\def\l@subsection{\@tocline{2}{0pt}{2.8pc}{5pc}{}}


\section{Introduction}

In this paper, we consider a simple dynamical system consisting of interacting particles moving on the sphere. Fix $d\geq2$, let $\S\subset\R^d$ be the unit sphere. 
The states of $n$ particles $z_1(t),\ldots,z_n(t)\in\S$ evolve according to
\begin{align}\label{eq:particle-system-eqn}
    \partial_tz_i(t)=\proj_{z_i(t)}\cM_t^n,
    \qquad
    \cM_t^n:=\frac1n\sum_{j=1}^nz_j(t) \,, 
    \qquad
    t \geq 0 \,, i=1,\ldots,n \, ,
\end{align}
where $\proj_x:=\Id-xx^\top$, so that $\proj_xy$ is the
orthogonal projection of $y\in\R^d$ onto $T_x\S$ for $x \in \S$.
These dynamics describe alignment with the empirical mean, with the projection ensuring that each particle remains on the sphere.

Equation~\eqref{eq:particle-system-eqn} can be regarded as a high-dimensional generalization of the classical Kuramoto model~\cite{kuramoto1975self}. Since its introduction in the 1970s, the Kuramoto model has served as a paradigmatic example of synchronization phenomena, with applications across physics~\cite{wiesenfeld1996synchronization,kozyreff2000global} and biology~\cite{goldstein2009noise,myung2015gaba}.
When $d=2$, equation~\eqref{eq:particle-system-eqn} reduces to the well-known identical-frequency Kuramoto model~\cite{acebron2005kuramoto,benedetto2015complete}.
When $d\geq 3$, the model captures the essential dynamics of orientation alignment in swarming and flocking~\cite{chandra2019continuous,frouvelle2019long}.
More recently, this model has also provided a tractable framework for understanding clustering and oversmoothing in transformers, a class of neural networks used in language modeling. 
In this setting, the alignment of interacting particles helps explain how data representations form clusters or lose their distinctiveness as they evolve across network layers~\cite{geshkovski2023emergence,geshkovski2023mathematical,criscitiello2024synchronization,bruno2024emergence,burger2025analysis,chen2025quantitative,andrew25,bruno2025multiscale,chen2026propagation}.

A central question in the study of this model concerns the long-time behavior of the interacting particles, including whether they synchronize, form clusters or align.
These phenomena are often analyzed through the corresponding mean-field equation, obtained in the limit as the number of particles tends to infinity~\cite{carrillo2014contractivity,morales2022trend}.
However, transferring conclusions about the long-time dynamics of this mean-field equation to finite particle systems calls for quantitative approximation estimates that remain uniform in time.
Establishing such estimates remains a significant analytical challenge, as it requires controlling the accumulation of errors over arbitrarily long times.

Our aim in this paper is therefore to quantify this difference uniformly over all time, a problem known as uniform-in-time propagation of chaos~\cite{sznitman1991topics}.
Let $\cP(\S)$ denote the space of Borel probability measures on $\S$, and define the empirical measure of the particle system by
\begin{align*}
\mu_t^n := \frac{1}{n}\sum_{i=1}^n \delta_{z_i(t)}.
\end{align*}
For independent initial states with common law $\mu_0\in\cP(\S)$, let $\mu_t\in\cP(\S)$ denote the corresponding mean-field evolution, defined below in~\eqref{eq:intro mean field}.
Under suitable assumptions on $\mu_0$, we seek bounds of the form
\begin{align*}
\mathbb P\!\left(
\sup_{t\geq0}W_2(\mu_t^n,\mu_t)\leq a_{n,\eta}
\right)\geq1-\eta,
\quad \text{with } a_{n,\eta}\longrightarrow0\quad\text{as }n\to\infty,
\end{align*}
for each fixed confidence parameter $\eta\in(0,1)$, where $W_2$ is the quadratic Wasserstein distance associated with the Euclidean metric restricted to $\S$.
Such an estimate controls the entire empirical-measure trajectory on a single event of high probability.
Achieving this estimate is more demanding than establishing probability bounds separately at each time, even when those bounds are uniform in time.

The mean-field equation of~\eqref{eq:particle-system-eqn} can be derived by writing the evolution equation for the empirical measure $\mu_t^n$, and formally taking the limit as $n\to\infty$.
Let $\sigma$ denote surface measure on $\S$. For $\mu\in\cP(\S)$, define
\begin{align}\label{eq:mean vector}
    \cM_\mu:=\int_{\S}y\,\diff\mu(y),
    \qquad
    \cX_\mu(x):=\proj_x\cM_\mu.
\end{align}
Then the mean-field evolution $\mu_t$ satisfies the continuity equation
\begin{align}\label{eq:intro mean field}
    \partial_t\mu_t+\dive(\mu_t\cX_{\mu_t})=0,
    \qquad \mu_{t=0}=\mu_0,
\end{align}
where $\dive=\mathrm{div}_{\S}$ denotes the divergence on the sphere.
The empirical measure $\mu_t^n$ satisfies the same equation in the weak sense, with initial measure $\mu_0^n$.
Thus, in our setting, uniform-in-time propagation of chaos can be viewed as a quantitative stability problem for the mean-field equation under empirical perturbations of the initial measure.

\textbf{Uniform-in-time propagation of chaos.}
Finite-time propagation of chaos can be obtained by standard stability estimates for the regular velocity field $\cX_\mu$, but such stability estimates usually deteriorate exponentially over time~\cite{dobrushin1979vlasov}.
To obtain uniform-in-time estimates, we need additional information of the dynamical system.
For stochastic systems, uniform-in-time estimates can be obtained by exploiting both the regularization and dissipation provided by diffusion~\cite{malrieu2001logarithmic,lacker2023sharp,delgadino2023phase,delarue2025uniform,rosenzweig2023global,gerber2025consensus}.
Without diffusion, such estimates can be obtained through stability or relaxation properties of the deterministic interacting particle system.
One approach uses contraction in transport distances, as in aggregation equations with convex interaction potentials~\cite{carrillo2011global} and the Kuramoto model confined to a sufficiently short arc~\cite{carrillo2014contractivity}.
More generally, uniform stability with respect to initial data yields results for Cucker--Smale systems in a flocking regime~\cite{ha2018uniform,gerber2025uniform} and a singular nonlocal transport model~\cite{agazzi2026quantitative}.
Other results exploit more specific structures: closed equations for the mean and covariance yield uniform Wasserstein estimates for Gaussian data in Stein variational gradient descent~\cite{liu2023towards}, while decay of the mean-field training loss, together with the kernel structure of the network, yields uniform weak estimates in neural-network training~\cite{glasgow2026uniform}.

The Kuramoto model considered here allows initial densities supported on the entire sphere.
This support persists at every finite time, so evolving the density does not bring the solution into the short-arc regime where contraction is available.
Moreover, stability estimates cannot be uniform both in time and over all admissible initial data, since arbitrarily small perturbations of the uniform distribution can synchronize in opposite directions. 
Finally, even when both systems synchronize, their limiting directions may differ, so synchronization alone does not yield the desired uniform-in-time estimate.

\textbf{Contributions.}
We establish uniform-in-time propagation of chaos in the quadratic Wasserstein distance for the Kuramoto model~\eqref{eq:particle-system-eqn} in every dimension $d\geq2$, assuming independent initialization from a law with an $L^p$ density, $p>1$, and a nonzero mean. In sufficiently high dimensions, our bound matches the empirical Wasserstein sampling rate $n^{-1/(d-1)}$ at each fixed confidence level and is therefore optimal in its power of $n$. To our knowledge, this is the first such result to allow initial densities supported on the entire sphere, beyond the short-arc regime covered by the contractivity estimates of~\cite{carrillo2014contractivity}.
Our proof builds on and extends the quantitative mean-field synchronization estimates~\cite{morales2022trend,chen2025quantitative}.
We also establish a late-time comparison estimate for the mean vectors of the particle and mean-field systems, allowing us to combine finite-time mean-field bounds with synchronization to obtain uniform-in-time control.

\textbf{Organization.}
\Cref{sec:main-results} states the main result \Cref{thm:uniform-time-poc-final} and its high-dimensional consequence \Cref{cor:empirical-rate-high-dimension}, and explains the proof strategy. 
We establish quantitative synchronization in \Cref{sec:synchronization} and stability near
synchronization in \Cref{sec:stability}. 
\Cref{sec:propagation-of-chaos} then develops the finite-time propagation of chaos and completes the proof of uniform-in-time propagation of chaos. 
The sampling and concentration estimates are collected in \Cref{app:sampling-estimates}.

\section{Main results}\label{sec:main-results}

Throughout the paper, $W_2$ denotes the quadratic Wasserstein distance for
the Euclidean metric restricted to $\S$:
\begin{align}\label{eq:def-wasserstein}
    W_2^2(\mu,\nu):=
    \inf_{\pi\in\Pi(\mu,\nu)}
    \int_{\S\times\S}\|x-y\|_2^2\,\diff\pi(x,y),
\end{align}
where $\Pi(\mu,\nu)\subset\mathcal P(\S\times\S)$ is the set of couplings of $\mu$ and $\nu$.
In particular, $W_2\leq2$. 

\subsection{Assumptions}
Let $(\mu_t)_{t\ge0}$ be the Kuramoto flow
with initial measure $\mu_0\in\cP(\S)$, that is, the weak solution
of \eqref{eq:intro mean field}.

\begin{assumption}\label{a:L^p mu_0}
The initial positions $z_1(0),\ldots,z_n(0)$ are i.i.d. with common law $\mu_0\in\cP(\S)$, where
$\diff\mu_0=f_0\,\diff\sigma$ and
$f_0\in L^p(\S,\diff\sigma)$ for some $p\in(1,\infty)$.
\end{assumption}

\begin{assumption}\label{a:R_0 > 0}
The initial mean vector satisfies $R_0\coloneqq \|\cM_{\mu_0}\|_2>0$.
\end{assumption}

For the $p$ in \Cref{a:L^p mu_0} and any $\varepsilon \in (0,1/2)$, set
\begin{align}\label{eq:def-beta-p-epsilon}
    p^*:=\frac{p}{p-1},\qquad
    \gamma_p:=\frac{d-1}{p^*},\qquad
    \beta_p:=\min\{2,\gamma_p\},\qquad
    \beta_{p,\varepsilon}:=(1-\varepsilon)\beta_p.
\end{align}
We also define the compact-space empirical $W_2^2$ rate~\cite{fournier2015rate,weed2019sharp}
\begin{align}\label{eq:def-empirical-rate}
    \mathfrak r_m(n):=
    \begin{cases}
        n^{-1/2},&m<4,\\
        n^{-1/2}\log(2+n),&m=4,\\
        n^{-2/m},&m>4.
    \end{cases}
\end{align}
In the following, the constants $C>0$ may depend on the full initial law $\mu_0$, as well as on $d$, $p$, and $\varepsilon$, but are independent of $n$, $\eta$, and time $t$.

\subsection{Uniform-in-time propagation of chaos}

\begin{theorem}\label{thm:uniform-time-poc-final}
Adopt \Cref{a:L^p mu_0,a:R_0 > 0}. For every
$\varepsilon\in(0,1/2)$, there exists
$C=C(\varepsilon,\mu_0,p,d)>0$ such that, for every
$n\geq2$ and $\eta\in(0,1/2)$, with probability at least
$1-\eta$,
\begin{align}\label{eq:uniform-time-poc-final}
    \sup_{t\geq0}W_2^2(\mu_t^n,\mu_t)
    \leq C\left[
        \mathfrak r_{d-1}(n)
        +\bigl(1+\log(8/\eta)\bigr)^{\alpha_{p,\varepsilon}}
          n^{-\alpha_{p,\varepsilon}}
    \right],
\end{align}
where
\begin{align}\label{eq:def-alpha-p-epsilon}
    \alpha_{p,\varepsilon}
    :=\frac{\beta_{p,\varepsilon}}{\beta_{p,\varepsilon}+2}
    <\frac12.
\end{align}
\end{theorem}

The bound in~\eqref{eq:uniform-time-poc-final} separates the empirical sampling contribution $\mathfrak r_{d-1}(n)$ from an additional term arising from the discrepancy in the synchronization direction.
The following corollary identifies a dimension range in which the empirical sampling term $\mathfrak r_{d-1}(n)$ dominates the rate in $n$ at fixed confidence.

\begin{corollary}\label{cor:empirical-rate-high-dimension}
Adopt \Cref{a:L^p mu_0,a:R_0 > 0}. Suppose that
\begin{align}\label{eq:high-dimensional-condition}
    d-1>2+\frac4{\beta_p}.
\end{align}
Then $\varepsilon\in(0,1/2)$ can be chosen, depending only on
$(p,d)$, so that $\alpha_{p,\varepsilon}>2/(d-1)$.
Consequently, there exists $C=C(\varepsilon,\mu_0,p,d)>0$
such that, for every integer $n\geq2$ and
$\eta\in(0,1/2)$, with probability at least $1-\eta$,
\begin{align}\label{eq:empirical-rate-high-dimension-w2}
    \sup_{t\geq0}W_2^2(\mu_t^n,\mu_t)
    \leq C\left(\frac{1+\log(8/\eta)}{n}\right)^{2/(d-1)}.
\end{align}
In the regime $\beta_p=2$, condition
\eqref{eq:high-dimensional-condition} reduces to $d-1>4$.
\end{corollary}

The proof combines a finite-time mean-field estimate with a late-time estimate that applies after synchronization has begun.
To see how the two parts fit together, fix $T\geq0$ and apply the comparison estimate in \Cref{cor:late-time-w2} to obtain
\begin{equation}\label{eq:proof-strategy-reduction}
    \sup_{t\geq0}W_2^2(\mu_t^n,\mu_t)
    \leq C_0 \left(\sup_{0\leq t\leq T}W_2^2(\mu_t^n,\mu_t)+1-R_T\right).
\end{equation}
Here $C_0$ is a universal constant, and
$R_T\coloneqq \|\cM_{\mu_T}\|_2$ approaches $1$ as synchronization occurs.
The uniform-in-time estimate thus reduces to a stability estimate up to time $T$ and the rate of synchronization.
In \Cref{sec:synchronization}, we show that $1-R_t$ decays exponentially after a burn-in time.
In~\Cref{sec:stability}, we bound the difference of the mean direction in terms of $1-R_T$, which yields \eqref{eq:proof-strategy-reduction}.
Finally, choosing an appropriate $T$ in~\Cref{cor:late-time-w2} yields the claimed uniform-in-time bound.
The detailed argument appears in \Cref{subsec:main-theorem-proof}.

\Cref{cor:empirical-rate-high-dimension} identifies a high-dimensional regime in which our uniform-in-time estimate is optimal in its dependence on $n$. Indeed, the absolute continuity of the initial law implies the lower bound $W_2^2(\mu_0^n,\mu_0)\ge c_{\mu_0,d} \, n^{-2/(d-1)}$, which also bounds the uniform-in-time error from below since the supremum includes $t=0$; see also \cite[Proposition~6]{weed2019sharp}.

One possible refinement for \Cref{thm:uniform-time-poc-final} and \Cref{cor:empirical-rate-high-dimension} is to make the dependence of the constant $C$ on the initial spread explicit. 
When the initial distribution is concentrated at a single point, the empirical and mean-field distributions coincide at all times. 
This motivates asking whether $C$ can be shown to decrease as the initial distribution becomes more concentrated, while preserving the same convergence rates in $n$.

A natural extension of our results concerns transformer dynamics. Such an extension would clarify when the mean-field approximation of a finite collection of tokens remains accurate as the network becomes deeper. Synchronization has already been established for transformer dynamics with small parameters~\cite{chen2025quantitative}, suggesting that our strategy of combining finite-time approximation with synchronization may also apply in this setting.

As a qualitative complement to the quantitative theorems, the
late-time comparison \Cref{cor:late-time-w2} also gives uniform-in-time stability at any
synchronizing solution.

\begin{corollary}[Uniform stability at synchronizing solutions]
\label{cor:qualitative-poc}
Let $(\mu_t)_{t\ge0}$ and, for each $k\ge1$,
$(\nu_t^k)_{t\ge0}$ be Kuramoto flows, with initial
measures $\mu_0$ and $\nu_0^k$ in $\cP(\S)$.
Suppose that $R_t \coloneqq \|\cM_{\mu_t}\|_2\to1$ as $t\to\infty$.
If $\lim_{k\to\infty}W_2(\nu_0^k,\mu_0) = 0$, then
\begin{align*}
    \lim_{k \to \infty} \sup_{t\geq0}W_2(\nu_t^k,\mu_t) =0.
\end{align*}
In particular, if $\mu_0$ is atomless and $R_0>0$, then, for a single
infinite i.i.d. sequence of initial particles with law $\mu_0$,
\begin{align*}
    \lim_{n \to \infty} \sup_{t\geq0}W_2(\mu_t^n,\mu_t)=0
    \qquad\text{almost surely}.
\end{align*}
\end{corollary}

\begin{proof}
Transporting an optimal coupling of $\nu_0^k$ and $\mu_0$ along their
characteristic flows and applying \Cref{lem:one-sided-kuramoto} gives
$W_2^2(\nu_t^k,\mu_t)\leq e^{2t}W_2^2(\nu_0^k,\mu_0)$.
Thus \Cref{cor:late-time-w2} implies, for every $T\geq0$,
\begin{align*}
    \sup_{t\geq0}W_2^2(\nu_t^k,\mu_t)
    \leq C\left(e^{2T}W_2^2(\nu_0^k,\mu_0)+1-R_T\right).
\end{align*}
The first assertion follows from letting $k\to\infty$ first and then $T\to\infty$.
For the last assertion, \cite[Theorem~1]{frouvelle2019long} gives
$R_t\to1$, while the empirical measures $\mu_0^n$ converge to
$\mu_0$ in $W_2$ almost surely.
\end{proof}

\section{Quantitative synchronization}\label{sec:synchronization}

We first establish a quantitative decay estimate for the synchronization deficit $1-R_t$. Recall that $\S\subseteq\R^d$ is the unit sphere ($d\geq2$) and that, for $\mu\in\cP(\S)$, \eqref{eq:mean vector} defines
    \begin{align*}
        \cM_{\mu} \coloneqq \int_{\S} y \diff \mu(y),  \text{ and }  \cX_{\mu}(x) \coloneq \proj_x \cM_{\mu} = \int_{\S} \proj_x[ y]  \diff \mu(y), 
    \end{align*}
where $\proj_x[ y]:= y -\langle x,y\rangle\, x$ denotes the projection onto the hyperplane $T_x\S$ tangent to the sphere at $x \in \S$, $\langle x,y\rangle$ denotes the standard inner product in $\R^d$ and $\|x\|_2 \coloneqq \sqrt{ \langle x, x \rangle}$. We also denote
    \begin{align}
        \cS_{\mu} \coloneqq \int_{\S} y y^{\top} \diff \mu(y) \in \R^{d \times d}.
    \end{align}

The Kuramoto evolution of $\mu_t$ is governed by the continuity equation
    \begin{align}\label{eq:main kuramoto}
            \partial_t \mu_t + \dive\left(\mu_t  \mathcal{X}_{\mu_t}\right) = 0, \quad  t \geq 0,
    \end{align}
where here and throughout the paper $\dive=\mathrm{div}_{\S}$ denotes the divergence operator on the sphere. For $\cX_{\mu}$ in \eqref{eq:mean vector}, we frequently use the fact that for any $ x \in \S$, 
    \begin{align}\label{eq:div X_mu}
        \dive \cX_{\mu} (x) = - (d-1) \langle \cM_{\mu} , x \rangle ,
    \end{align}
which follows from a direct computation. For simplicity of notation, set
    \begin{align}
        \cS_t \coloneqq \cS_{\mu_t}, \quad  \cM_t \coloneqq \cM_{\mu_t}, \quad \cX_t \coloneqq \cX_{\mu_t}, \quad R_t \coloneq \|\cM_t\|_2, \quad \cU_t \coloneq \frac{\cM_t}{R_t} \in \S,
    \end{align}
where $\cU_t$ is well-defined whenever $R_t \neq 0$. Indeed, \Cref{lemma:derivative M_t R_t U_t} guarantees that $R_t$ is non-decreasing, and thus positive so long as $R_0 > 0$.

\begin{lemma}[{\cite[Lemma 5.2, Lemma 6.2]{chen2025quantitative}}]\label{lemma:derivative M_t R_t U_t}
    For the derivatives of $\cM_t$ and $R_t$, we have the following formulas:
        \begin{align*}
            \partial_t \cM_t = \int_{\S} \cX_t(y) \diff \mu_t(y) = \cM_t - \cS_t \cM_t,
        \end{align*}
    and 
        \begin{align*}
            \partial_t (R_t ^2) = 2 \int_{\mathbb{S}^{d-1}}  \| \cX_t(y)\|_2 ^2  \diff \mu_t(y),
        \end{align*}
    Thus, $R_t$ is non-decreasing. Moreover, if $R_0>0$, then
        \begin{align*}
            \partial_t \cU_t = \frac{1}{R_t} \proj_{\cU_t} \left[ \partial_t \cM_t \right] = - \proj_{\cU_t} [\cS_t \cU_t] .
        \end{align*}
    
\end{lemma}
Next, we define the spherical caps with centers $\pm U \in \S$ for $\alpha \in (0, \pi /2)$,
\begin{align}
    S_{\alpha} ^+ (U) \coloneq \left\{ x \in \mathbb{S}^{d-1} \ | \ \left\langle x , U \right\rangle \geq \cos \alpha \right\}, \quad S_{\alpha} ^- (U) \coloneq \left\{ x \in \mathbb{S}^{d-1} \ | \ \left\langle x , -U \right\rangle \geq \cos \alpha \right\} . 
\end{align}
The spherical caps become smaller as $\alpha \to 0$. 
    For simplicity, we also write $S_{\alpha} ^+(t) = S_{\alpha} ^+ (\cU_t)$ and $S_{\alpha} ^-(t) = S_{\alpha} ^- (\cU_t)$.
Our main result in this section establishes exponential decay of the
synchronization deficit $1-R_t$ at rates arbitrarily close to
$\beta_p=\min\{2,\gamma_p\}$, defined in
\eqref{eq:def-beta-p-epsilon}.

\begin{proposition}\label{prop:sharp exponential rate R_t}
Adopt \Cref{a:L^p mu_0,a:R_0 > 0}.
Fix $\varepsilon\in(0,1/2]$, and let $t_\varepsilon\geq0$
be chosen so that
\eqref{eq:exponential decay negative cap initial condition} holds.
Then, for any $t\geq t_\varepsilon$,
\begin{align}\label{eq:sharp exponential rate R_t derivative}
    \partial_t R_t
    \geq 2(1-\varepsilon)(1-R_t)
    -4\mu_t\left(\S\setminus S_\varepsilon^+(t)\right).
\end{align}
Moreover, with $T_\varepsilon$ defined in \eqref{eq: T_eps},
there exists
$C_\varepsilon=C_\varepsilon(d,\mu_0,p,\varepsilon)>0$
such that, for any $t\geq t_\varepsilon+T_\varepsilon$,
\begin{align}\label{eq:sharp exponential rate R_t}
    1-R_t
    \leq
    \left(1+C_\varepsilon(t-t_\varepsilon-T_\varepsilon)\right)
    e^{-(1-\varepsilon)\beta_p
       (t-t_\varepsilon-T_\varepsilon)}.
\end{align}
\end{proposition}

The proof begins with the following qualitative synchronization result.

\begin{proposition}[{\cite[Theorem~1]{frouvelle2019long}}]
\label{prop:kuramoto synchronize}
Adopt \Cref{a:L^p mu_0,a:R_0 > 0}.
Then the solution $\mu_t$ of \eqref{eq:main kuramoto}
converges weakly, as $t\to\infty$, to a Dirac mass
$\delta_{\cU_\infty}$ for some $\cU_\infty\in\S$.
\end{proposition}

\Cref{lemma:derivative M_t R_t U_t} and \Cref{prop:kuramoto synchronize} imply
\begin{align}\label{eq:vanishing derivatives Rt Ut}
    R_t\to 1
    \qquad\text{as }t\to\infty.
\end{align}
Consequently, for every $\varepsilon\in(0,1/2]$, there exists
$t_\varepsilon\geq0$ such that
$R_t\geq1-\varepsilon/2$ for all $t\geq t_\varepsilon$,
as required in
\eqref{eq:exponential decay negative cap initial condition}.

For initial densities $f_0$ in $L^2$,
\cite[Theorem~2.5]{chen2025quantitative} already provides
a quantitative exponential synchronization estimate for $1-R_t$ with a rate $1/100$.
For our argument in this paper, however, qualitative convergence suffices
to guarantee the existence of $t_\varepsilon$ as required in
\eqref{eq:exponential decay negative cap initial condition}.
Starting from this time, we derive quantitative estimates
that yield the decay bound in
\Cref{prop:sharp exponential rate R_t}.
To this end, we study the characteristic flow generated by
$\cX_t$ and the mass near the antipodal direction. For $s\geq0$, let $\phi_{s\to t}:\S\to\S$ be the
characteristic flow defined by
    \begin{align}\label{eq:C1 diffeomorphism}
        \partial_t \phi_{s\to t}(x) = \cX_{t}(\phi_{s\to t}(x)), \text{ with }  \phi_{s \to s}(x) = x, \ \forall x \in \S .
    \end{align}
Uniqueness of the ODE flow implies that, for any $s,r,t\geq0$,
    \begin{align*}
        \phi_{r \to t} \circ \phi_{s \to r} = \phi_{s \to t}.
    \end{align*}
In particular, for any $s,t \geq 0$, $\phi_{s \to t} \circ \phi_{t \to s} = \Id$.

By the characteristic-flow formula in
\cite[proof of Theorem 1.3, (2.3)--(2.4)]{lin2026sharp},
we may write $\diff\mu_t=f_t\,\diff\sigma$ with
$f_t\in L^p(\S,\diff\sigma)$ for every finite $t$.
For any Borel set $A\subseteq\S$, define
    \begin{align}
        f_t ^p (A) \coloneqq \int_A f_t ^p(x) \diff \sigma(x).
    \end{align}
Recall from \eqref{eq:def-beta-p-epsilon} that $p^{\ast}=p/(p-1)$,
so $1/p+1/p^{\ast}=1$ and
    \begin{align*}
        \gamma_p = (d-1)/p^{\ast} .
    \end{align*}
The following estimate controls the motion of the mean direction by the synchronization deficit.
\begin{lemma}\label{lem: bound derivative Ut}
    Assume $R_0>0$. Then, for any $t\geq0$,
        \begin{align}\label{eq:Ft bound derivative Ut}
         \|\partial_t \cU_t\|_2 \leq  (1-R_t).
        \end{align}
\end{lemma}
\begin{proof}
    By \Cref{lemma:derivative M_t R_t U_t}, we can write
        \begin{align*}
            \partial_t \cU_t = - \proj_{\cU_t} [\cS_t \cU_t]  =   -\int_{\S} \langle y ,\cU_t \rangle \proj_{\cU_t} y  \diff \mu_t(y) = \int_{\S} \left(1-\langle y ,\cU_t \rangle \right)\proj_{\cU_t} y  \diff \mu_t(y),
        \end{align*}
    where we used the fact that $\int_{\S} \proj_{\cU_t} y  \diff \mu_t(y) = \proj_{\cU_t} \cM_t = 0$. \eqref{eq:Ft bound derivative Ut} then follows from $\|\proj_{\cU_t} y\|_2 \leq 1$. 
        \begin{align*}
            \|\partial_t \cU_t\|_2 \leq  \int_{\S} \left(1-\langle y ,\cU_t \rangle \right)  \diff \mu_t(y) = 1- R_t.
        \end{align*}
\end{proof}

\begin{proposition}\label{prop:exponential decay negative cap}
    Adopt \Cref{a:L^p mu_0} and \Cref{a:R_0 > 0}. Then for any $\varepsilon \in (0,1/2]$, there is a $t_{\varepsilon} = t_{\varepsilon} (d,\mu_0,p, \varepsilon)$ and a $C = C(d,\mu_0,p)$, such that for any $t \geq t_{\varepsilon}$,
        \begin{align}\label{eq:f^p decay negative cap}
            f_t ^p \left(S_{\varepsilon} ^-(t) \right) \leq C \cdot e^{-(t-t_{\varepsilon}) \cdot (1-\varepsilon)p \gamma_p },
        \end{align}
    and for any $\delta\in(0,\varepsilon]$ and $t\geq t_\varepsilon$, 
        \begin{align}\label{eq:mu_t decay negative cap}
            \mu_t \left(S_{\delta} ^-(t) \right) \leq C \cdot \delta^{\gamma_p } \cdot e^{-(t-t_{\varepsilon})\cdot (1-\varepsilon) \gamma_p}.
        \end{align}
\end{proposition}
\begin{proof}
    Since $R_t\to1$ by \eqref{eq:vanishing derivatives Rt Ut}, choose
$t_\varepsilon=t_\varepsilon(d,\mu_0,p,\varepsilon)$ so that when $t \geq t_{\varepsilon}$,
\begin{align}\label{eq:exponential decay negative cap initial condition}
    R_t\geq1-\varepsilon/2.
\end{align}
Fix $\varepsilon_1\in(\varepsilon,\pi/2)$ and choose a smooth
nondecreasing function
$\xi_{\varepsilon,\varepsilon_1}:\R\to[0,1]$ such that
$\xi_{\varepsilon,\varepsilon_1}(a)=0$ for
$a\leq\cos\varepsilon_1$ and
$\xi_{\varepsilon,\varepsilon_1}(a)=1$ for
$a\geq\cos\varepsilon$. We may choose it with
$0\leq\xi_{\varepsilon,\varepsilon_1}'(a)
\leq2/(\cos\varepsilon-\cos\varepsilon_1)$.
We will let $\varepsilon_1 \to \varepsilon$ at the end.

    Differentiating along characteristics yields, for almost every
$t\geq t_\varepsilon$,
        \begin{align}\label{eq:derivative f^p negative cap}
            \begin{split}
                &\frac{\diff}{\diff t} \int_{\S} \xi_{\varepsilon,\varepsilon_1} \left( \langle x, -\cU_t \rangle \right) f_t ^p (x) \diff \sigma(x) 
                \\ &= \int_{\S} \xi_{\varepsilon,\varepsilon_1} ' \left( \langle x, -\cU_t \rangle \right)  \langle x , - \partial_t \cU_t \rangle f_t ^p (x) +  \xi_{\varepsilon,\varepsilon_1} \left( \langle x, -\cU_t \rangle \right) p f_t ^{p-1}(x) \partial_t f_t(x) \diff \sigma(x)
                \\ &= \int_{\S} \xi_{\varepsilon,\varepsilon_1} ' \left( \langle x, -\cU_t \rangle \right)  \left[\langle x , - \partial_t \cU_t \rangle - R_t \|\proj_x \cU_t \|_2 ^2\right] f_t ^p (x) 
                \\  & \quad +  \xi_{\varepsilon,\varepsilon_1} \left( \langle x, -\cU_t \rangle \right) (p-1)(d-1) R_t\langle x , \cU_t \rangle f_t ^p (x) \diff \sigma(x).
            \end{split}
        \end{align}
    The factor $\xi_{\varepsilon,\varepsilon_1}'(\langle x,-\cU_t\rangle)$
is supported where
$\cos\varepsilon_1\leq\langle x,-\cU_t\rangle\leq\cos\varepsilon$.
On this set, $\langle\cU_t,\partial_t\cU_t\rangle=0$ and
$\|\proj_{\cU_t}x\|_2\leq\sin\varepsilon_1$, so
        \begin{align*}
            \langle x , - \partial_t \cU_t \rangle = \langle \proj_{\cU_t} x , - \partial_t \cU_t \rangle \leq \sin(\varepsilon_1) \| \partial_t \cU_t \|_2 \leq (1-R_t) \sin(\varepsilon_1)  \leq \frac{\varepsilon}{2} \sin(\varepsilon_1)
        \end{align*}
    and
        \begin{align*}
            -R_t \|\proj_x \cU_t \|_2 ^2 \leq -R_t \sin^2(\varepsilon) \leq -(1-\varepsilon/2)\sin^2(\varepsilon).
        \end{align*}
    Since $\varepsilon\leq1/2$, choosing $\varepsilon_1>\varepsilon$ sufficiently close to $\varepsilon$ makes the first integral on the right-hand side of \eqref{eq:derivative f^p negative cap} nonpositive. On the support of
$\xi_{\varepsilon,\varepsilon_1}(\langle x,-\cU_t\rangle)$,
        \begin{align*}
            R_t\langle x , \cU_t \rangle \leq -(1-\varepsilon/2)\cos(\varepsilon_1) \leq -(1-\varepsilon)
        \end{align*}
    for $\varepsilon_1$ sufficiently close to $\varepsilon$. Hence, we obtain that when $t \geq t_{\varepsilon}$,
        \begin{align*}
            &\frac{\diff}{\diff t} \int_{\S} \xi_{\varepsilon,\varepsilon_1} \left( \langle x, -\cU_t \rangle \right) f_t ^p (x) \diff \sigma(x) \\
            &\qquad\leq -(1-\varepsilon)(p-1)(d-1) \int_{\S} \xi_{\varepsilon,\varepsilon_1} \left( \langle x, -\cU_t \rangle \right) f_t ^p (x) \diff \sigma(x).
        \end{align*}
    Integrating from $t_{\varepsilon}$ to $t$, a Gr{\"o}nwall argument shows that
        \begin{align*}
            &\int_{\S} \xi_{\varepsilon,\varepsilon_1} \left( \langle x, -\cU_t \rangle \right) f_t ^p (x) \diff \sigma(x) \\
            &\qquad\leq e^{-(1-\varepsilon)(p-1)(d-1)(t-t_{\varepsilon})} \int_{\S} \xi_{\varepsilon,\varepsilon_1} \left( \langle x, -\cU_{t_{\varepsilon}} \rangle \right) f_{t_{\varepsilon}} ^p (x) \diff \sigma(x).
        \end{align*}
    Notice that $(p-1)(d-1) = p \gamma_p$ by definition. We then let $\varepsilon_1 \to \varepsilon$ on both sides, and see that for any $t \geq t_{\varepsilon}$,
        \begin{align}\label{eq:f^p decay negative cap 1}
            f_t ^p \left(S_{\varepsilon} ^-(t) \right) \leq e^{-(1-\varepsilon)p \gamma_p (t-t_{\varepsilon})} f_{t_{\varepsilon}} ^p \left(S_{\varepsilon} ^-(t_{\varepsilon}) \right).
        \end{align}

    To proceed, we first take $\varepsilon=1/2$. By \eqref{eq:f^p decay negative cap 1}, we can choose $t_{1/2}=t_{1/2}(d,\mu_0,p)$ so that, for $t\geq t_{1/2}$,
        \begin{align}
            f_t ^p \left(S_{1/2} ^-(t) \right) \leq e^{-(1/2)p \gamma_p (t-t_{1/2})} f_{t_{1/2}} ^p \left(S_{1/2} ^-(t_{1/2}) \right).
        \end{align}
        The characteristic-flow formula in \cite[proof of Theorem 1.3, (2.3)--(2.4)]{lin2026sharp} bounds
    $\|f_t\|_{L^p}$ on $[0,t_{1/2}]$. Together with the preceding
    decay estimate, this gives
        \begin{align*}
            f_t^p(S_{1/2}^-(t))\leq C(d,\mu_0,p) ,
        \end{align*}
    for all $t\geq0$.
    Since $S_\varepsilon^-(t_\varepsilon)
    \subseteq S_{1/2}^-(t_\varepsilon)$,
    \eqref{eq:f^p decay negative cap 1} proves
    \eqref{eq:f^p decay negative cap}.

    Finally, for \eqref{eq:mu_t decay negative cap}, we combine \eqref{eq:f^p decay negative cap} and H{\"o}lder's inequality, and obtain that
        \begin{align*}
            \mu_t \left(S_{\delta} ^-(t) \right) \leq \left[f_t ^p \left(S_{\varepsilon} ^-(t) \right) \right]^{\frac{1}{p}} \cdot \left[\sigma \left(S_{\delta} ^-(t) \right) \right]^{\frac{1}{p^{\ast}}}
                \leq C  \cdot e^{-(t-t_{\varepsilon})\cdot (1-\varepsilon) \gamma_p} \cdot \delta^{\gamma_p},
        \end{align*}
    for some $C= C(d,\mu_0 ,p) $. Here, we used the fact that $S^- _{\delta}(t) \subseteq S^- _{\varepsilon}(t) $ and $\sigma \left( S^- _{\delta}(t) \right)$ is the area of a geodesic ball in $\S$ of radius $\delta$ and hence of order $\delta^{d-1}$.

\end{proof}

To turn the antipodal-cap estimate into decay of $1-R_t$, we show that
once $1-R_t=O(\varepsilon)$, characteristics outside $S_\varepsilon^-(t)$
reach $S_\varepsilon^+(t)$ in $O(\log(1/\varepsilon))$ time and remain there.

\begin{proposition}\label{prop: admit non admit capture time}
Fix $\varepsilon\in(0,1/2]$. Let $(\mu_t)_{t\geq0}$ be a Kuramoto
flow \eqref{eq:main kuramoto} such that, for some $t_\varepsilon\geq0$,
\begin{align*}
    R_t\geq1-\frac{\varepsilon}{2},
    \quad \text{when } t\geq t_\varepsilon.
\end{align*}
Define
\begin{align}\label{eq: T_eps}
    T_\varepsilon:=10\log(2/\varepsilon).
\end{align}
Then, for every $t\geq t_\varepsilon$,
\begin{align*}
    \phi_{t\to t+T_\varepsilon}
    \bigl(\S\setminus S_\varepsilon^-(t)\bigr)
    \subseteq
    S_\varepsilon^+(t+T_\varepsilon).
\end{align*}
Consequently,
\begin{align*}
    \mu_{t+T_\varepsilon}
    \bigl(\S\setminus S_\varepsilon^+(t+T_\varepsilon)\bigr)
    \leq
    \mu_t\bigl(S_\varepsilon^-(t)\bigr).
\end{align*}
Moreover, for every $s\geq t\geq t_\varepsilon$,
\begin{align*}
    \phi_{t\to s}\bigl(S_\varepsilon^+(t)\bigr)
    \subseteq S_\varepsilon^+(s),\qquad 
    \phi_{t\to s}\bigl(\S\setminus S_\varepsilon^-(t)\bigr)
    \subseteq\S\setminus S_\varepsilon^-(s).
\end{align*}
\end{proposition}

\begin{proof}
Fix $t\geq t_\varepsilon$ and set
\begin{align*}
    X_r:=\phi_{t\to r}(x),
    \qquad
    h_r:=\langle X_r,\cU_r\rangle,
    \qquad r\geq t.
\end{align*}
Whenever
\begin{align*}
    -\cos\varepsilon\leq h_r\leq\cos\varepsilon,
\end{align*}
set $s_r:=\sqrt{1-h_r^2}$. Then $s_r\geq\sin\varepsilon$, and
\Cref{lem: bound derivative Ut} gives
\begin{align*}
    \partial_rh_r =
    R_r(1-h_r^2)
    +\langle\proj_{\cU_r}X_r,\partial_r\cU_r\rangle \geq
    R_rs_r^2-s_r(1-R_r)
    =
    s_r\bigl(R_rs_r-(1-R_r)\bigr).
\end{align*}
Since $\varepsilon\leq1/2$,
\begin{align*}
    \sin\varepsilon
    \geq\varepsilon-\frac{\varepsilon^3}{6}
    \geq\frac{23}{24}\varepsilon.
\end{align*}
Since $\varepsilon/2\leq(12/23)s_r$ and $R_r\geq3/4$,
\begin{align}\label{eq:middle-region-drift}
    \partial_rh_r
    \geq s_r\left(\frac34s_r-\frac\varepsilon2\right)
    \geq\left(\frac34-\frac{12}{23}\right)s_r^2
    \geq\frac15(1-h_r^2).
\end{align}
At the boundary $h_r=\cos\varepsilon$,
\eqref{eq:middle-region-drift} points into
$S_\varepsilon^+(r)$, so this moving cap is forward invariant.
At $h_r=-\cos\varepsilon$, the same estimate shows that
$\S\setminus S_\varepsilon^-(r)$ is forward invariant.

Suppose now that $x\notin S_\varepsilon^-(t)$. If
$x\in S_\varepsilon^+(t)$, the conclusion follows from forward
invariance. Otherwise,
\begin{align*}
    -\cos\varepsilon<h_t<\cos\varepsilon.
\end{align*}
Until $X_r$ enters $S_\varepsilon^+(r)$,
\eqref{eq:middle-region-drift} implies
\begin{align*}
    \frac{\mathrm d}{\mathrm dr}\operatorname{arctanh}h_r
    \geq\frac15.
\end{align*}
The entrance time is therefore at most
\begin{align*}
    10\operatorname{arctanh}(\cos\varepsilon)
    =10\log\cot\frac\varepsilon2
    \leq10\log\frac2\varepsilon=T_\varepsilon.
\end{align*}
This proves the first inclusion.

Finally, $\mu_{t+T_\varepsilon}
    =
    (\phi_{t\to t+T_\varepsilon})_\#\mu_t$, which gives the measure estimate.
\end{proof}
By \Cref{prop: admit non admit capture time,prop:exponential decay negative cap}, the mass outside the forward cap decays exponentially. This yields the desired bound on $1-R_t$ in \Cref{prop:sharp exponential rate R_t}.

\begin{proof}[Proof of \Cref{prop:sharp exponential rate R_t}]
    By \Cref{lemma:derivative M_t R_t U_t} and $\| \proj_{y} \cU_t\|_2 ^2 = 1-\langle y , \cU_t \rangle ^2 $, we can write
        \begin{align*}
            \begin{split}
                \partial_t  R_t &= R_t\int_{\S}  \| \proj_{y} \cU_t\|_2 ^2  \diff \mu_t(y) \geq  R_t\int_{S_{\varepsilon} ^+ (t)} \left( 1-\langle y , \cU_t \rangle \right) \left( 1+\langle y , \cU_t \rangle \right)   \diff \mu_t(y)
                \\ &\geq R_t(1+\cos\varepsilon)
                    \int_{S_\varepsilon^+(t)}
                    (1-\langle y,\cU_t\rangle)\,\diff\mu_t(y)
                \\  & \geq R_t(1+\cos(\varepsilon))\int_{\S} \left( 1-\langle y , \cU_t \rangle \right)   \diff \mu_t(y) - 4\mu_t \left( \S \setminus  S_{\varepsilon} ^+ (t)\right)
                \\  & \geq (1 - \varepsilon/2)(1+\cos(\varepsilon))(1-R_t)- 4\mu_t \left( \S \setminus  S_{\varepsilon} ^+ (t)\right) 
                \\ & \geq 2 (1 - \varepsilon) (1 - R_t) - 4\mu_t \left( \S \setminus  S_{\varepsilon} ^+ (t)\right) ,
            \end{split}
        \end{align*}
    where the last two inequalities use $R_t \geq R_{t_\varepsilon} \geq 1-\varepsilon/2$ and $(1 - \varepsilon/2) (1 + \cos(\varepsilon)) \geq 2(1 - \varepsilon)$ for $\varepsilon \in (0, 1/2]$. 

    To proceed, we see that when $t \geq t_{\varepsilon} + T_{\varepsilon}$, according to \Cref{prop: admit non admit capture time} and \Cref{prop:exponential decay negative cap}, for some $C = C(d,\mu_0,p)$
        \begin{align*}
            \mu_t \left( \S \setminus  S_{\varepsilon} ^+ (t)\right) \leq  \mu_{t-T_{\varepsilon}}\left( S_{\varepsilon} ^- (t-T_{\varepsilon}) \right) \leq C \cdot \varepsilon^{\gamma_p } \cdot e^{-(t-t_{\varepsilon}-T_{\varepsilon})\cdot (1-\varepsilon) \gamma_p}.
        \end{align*}
    A Gr{\"o}nwall argument then yields \eqref{eq:sharp exponential rate R_t}.
\end{proof}

\section{Stability near synchronization}\label{sec:stability}
\label{sec:comparison-proofs}

The exponential decay for $1-R_t$ established in \Cref{prop:sharp exponential rate R_t} uses $L^p$ regularity and
therefore does not apply to empirical measures. Once $R_t$ is
sufficiently close to one, \Cref{prop: admit non admit capture time}
shows that particles away from the antipode join the main cluster
within a bounded time. Atoms near the antipode, however, may persist
indefinitely. For the late-time comparison, we need to control the remaining
motion of the mean direction without assuming that an atomic flow
synchronizes. The following lemma uses an angular dissipation
functional to bound this motion and the displacement of the mean
vector in terms of $1-R_T$.

\begin{lemma}\label{lem:finite-variation-center}
Let $(\nu_t)_{t\geq T}$ be a Kuramoto flow satisfying
\eqref{eq:main kuramoto}, and write
\begin{align*}
    M_t^\nu:=\cM_{\nu_t}=R_t^\nu\cU_t^\nu.
\end{align*}
If $R_T^\nu>1/2$, then
\begin{align}\label{eq:Ut-total-variation}
    \int_T^\infty \|\partial_t\cU_t^\nu\|_2\,\mathrm dt
    \leq
    \frac{\pi}{2(2R_T^\nu-1)}(1-R_T^\nu),
\end{align}
and
\begin{align}\label{eq:Mt-total-variation}
    \sup_{t\geq T}\|M_t^\nu-M_T^\nu\|_2
    \leq
    \left(1+\frac{\pi}{2(2R_T^\nu-1)}\right)(1-R_T^\nu).
\end{align}
\end{lemma}

\begin{proof}
Let
\begin{align*}
    \theta_t(x):=\arccos\langle x,\cU_t^\nu\rangle
\end{align*}
be the angle between $x$ and the mean direction $\cU_t^\nu$. When
$0<\theta_t(x)<\pi$, write
\begin{align*}
    \omega_t(x):=
    \frac{\proj_{\cU_t^\nu}x}{\sin\theta_t(x)},
\end{align*}
and set $\omega_t(x)=0$ at the two endpoints $\pm \cU_t^\nu$. Since
$\int_\S\proj_{\cU_t^\nu}x\,\nu_t(\mathrm dx)=0$,
\Cref{lemma:derivative M_t R_t U_t} gives
\begin{align*}
    \partial_t\cU_t^\nu
    =\int_\S
       (1-\cos\theta_t(x))\sin\theta_t(x)\,
       \omega_t(x)\,\nu_t(\mathrm dx).
\end{align*}
Thus, 
\begin{align*}
    \mathcal A_t:=
    \int_\S
       (1-\cos\theta_t(x))\sin\theta_t(x)\,
       \nu_t(\mathrm dx) \geq \|\partial_t\cU_t^\nu\|_2 .
\end{align*}
To prove \eqref{eq:Ut-total-variation}, we control the accumulated value of $\mathcal A_t$.
For the Kuramoto flow $(\nu_t)_{t\geq T}$, let $\phi^\nu_{T\to t}$ denote its
characteristic flow and write $X_t(x):=\phi^\nu_{T\to t}(x)$. When $0<\theta_t(X_t)<\pi$, we have
\begin{align*}
    & -\sin\theta_t(X_t)\frac{\mathrm d}{\mathrm dt}\theta_t(X_t)
    =\frac{\mathrm d}{\mathrm dt}
      \langle X_t,\cU_t^\nu\rangle=\langle\proj_{X_t}M_t^\nu,\cU_t^\nu\rangle
      +\langle X_t,\partial_t\cU_t^\nu\rangle\\
    &=R_t^\nu\sin^2\theta_t(X_t)
      +\sin\theta_t(X_t)
       \langle\omega_t(X_t),\partial_t\cU_t^\nu\rangle.
\end{align*}
Dividing by $\sin\theta_t(X_t)$ yields
\begin{align}\label{e:angle derivative}
    \frac{\mathrm d}{\mathrm dt}\theta_t(X_t)
    =-R_t^\nu\sin\theta_t(X_t)
     -\langle\omega_t(X_t),\partial_t\cU_t^\nu\rangle.
\end{align}
The first term pulls the particle toward the mean direction; the
second accounts for movement of the direction itself. To measure
their combined effect, set
\begin{align*}
    \Phi(\theta):=\theta-\sin\theta,
    \qquad
    \mathcal T_t:=
    \int_\S\Phi(\theta_t(x))\,\nu_t(\mathrm dx).
\end{align*}
Differentiating $\mathcal T_t$ along characteristics and using
\eqref{e:angle derivative}, we obtain, for almost every $t\geq T$,
\begin{align*}
    \partial_t\mathcal T_t
    &=\int_\S \Phi'(\theta_t(x))
      \left(
        -R_t^\nu\sin\theta_t(x)
        -\langle\omega_t(x),\partial_t\cU_t^\nu\rangle
      \right)\nu_t(\mathrm dx)\\
    &=-R_t^\nu\mathcal A_t
      -\left\langle\partial_t\cU_t^\nu,\,
        \int_\S(1-\cos\theta_t(x))
        \omega_t(x)\,\nu_t(\mathrm dx)
       \right\rangle,
\end{align*}
where the last equality uses $\Phi'(\theta)=1-\cos\theta$ and the
definition of $\mathcal A_t$. Since $\|\omega_t(x)\|_2\leq1$ and
$M_t^\nu=R_t^\nu\cU_t^\nu$,
\begin{align*}
    \left\|
      \int_\S(1-\cos\theta_t(x))
      \omega_t(x)\,\nu_t(\mathrm dx)
    \right\|_2
    &\leq\int_\S(1-\cos\theta_t(x))\,\nu_t(\mathrm dx)\\
    &=1-\langle M_t^\nu,\cU_t^\nu\rangle
     =1-R_t^\nu.
\end{align*}
Consequently, for almost every $t\geq T$,
\begin{align*}
    \partial_t\mathcal T_t
    \leq -R_t^\nu\mathcal A_t
          +(1-R_t^\nu)\|\partial_t\cU_t^\nu\|_2
    \leq -(2R_t^\nu-1)\mathcal A_t
    \leq -(2R_t^\nu-1)\|\partial_t\cU_t^\nu\|_2.
\end{align*}
So, for every $S\geq T$,
\begin{align}\label{e:integral derivative U_t bound}
    (2R_T^\nu-1)
    \int_T^S\|\partial_t\cU_t^\nu\|_2\,\mathrm dt
    \leq\mathcal T_T-\mathcal T_S
    \leq\mathcal T_T.
\end{align}
Now for any $\theta \in [0,\pi]$,
\begin{align*}
    \Phi(\theta)\leq\frac{\theta}{2}(1-\cos\theta) \leq \frac{\pi}{2}(1-\cos\theta),
\end{align*}
where the first inequality follows from
\begin{align*}
\frac{\theta}{2}(1-\cos\theta)-(\theta-\sin\theta)
=\sin\theta-\frac{\theta}{2}(1+\cos\theta) =2\cos\frac{\theta}{2}
\left(\sin\frac{\theta}{2}-\frac{\theta}{2}\cos\frac{\theta}{2}\right) \geq 0.
\end{align*}
It follows that
\begin{align*}
    \mathcal T_T
    \leq\frac{\pi}{2}
       \int_\S(1-\cos\theta_T(x))\,\nu_T(\mathrm dx)
    =\frac{\pi}{2}(1-R_T^\nu).
\end{align*}
Combining this with \eqref{e:integral derivative U_t bound} and letting $S\to\infty$, we obtain \eqref{eq:Ut-total-variation}. 

Finally, since $M_t^\nu=R_t^\nu\cU_t^\nu$, we have
\begin{align*}
    M_S^\nu-M_T^\nu
    =(R_S^\nu-R_T^\nu)\cU_T^\nu
      +R_S^\nu(\cU_S^\nu-\cU_T^\nu).
\end{align*}
Since $R_t^\nu$ is nondecreasing and $R_S^\nu\leq1$, it follows that
\begin{align*}
    \begin{split}
        \|M_S^\nu-M_T^\nu\|_2
    &\leq R_S^\nu-R_T^\nu
       +R_S^\nu\|\cU_S^\nu-\cU_T^\nu\|_2\leq R_S^\nu-R_T^\nu
       +\int_T^S\|\partial_t\cU_t^\nu\|_2\,\mathrm dt 
       \\ &\leq 1-R_T^\nu
       +\int_T^\infty\|\partial_t\cU_t^\nu\|_2\,\mathrm dt.
    \end{split}
\end{align*}
This gives \eqref{eq:Mt-total-variation}.
\end{proof}

We use this finite-variation bound to compare two flows in terms of their distance at time $T$ and the synchronization deficit of only one flow.

\begin{proposition}\label{cor:late-time-w2}
For any two Kuramoto flows $(\mu_t)_{t\geq0}$ and
$(\nu_t)_{t\geq0}$ satisfying \eqref{eq:main kuramoto}, and any
$T\geq0$, there is a universal constant $C$ such that
\begin{align*}
    \sup_{t\geq T}W_2^2(\mu_t,\nu_t)
    \leq C\bigl(W_2^2(\mu_T,\nu_T)+1-R_T^\mu\bigr).
\end{align*}
\end{proposition}

\begin{proof}
Set $Q:=W_2^2(\mu_T,\nu_T)$,
$a_\mu:=1-R_T^\mu$, and $a_\nu:=1-R_T^\nu$.
If $Q+a_\mu>1/8$, the conclusion follows from $W_2\leq2$
on $\S$. We may therefore assume $Q+a_\mu\leq1/8$.
Using the product coupling and the monotonicity of both $R_t^\mu$ and $R_t^\nu$ in $t$,
we obtain, for $t\geq T$,

\begin{equation}\label{eq: W2 bound}
\begin{split}
    W_2^2(\mu_t, \nu_t) & \leq 1 - (R_t^\mu)^2 + 1 - (R_t^\nu)^2 + \|M_t^\mu - M_t^\nu\|_2^2 \\
    & \leq 1 - (R_T^\mu)^2 + 1 - (R_T^\nu)^2 + \|M_t^\mu - M_t^\nu\|_2^2 \\
    & \leq 2 (a_\mu + a_\nu) + \|M_t^\mu - M_t^\nu\|_2^2 
\end{split}
\end{equation}
To bound the first term on the right hand side, note that for any
$x\in\S$ and $\rho\in\cP(\S)$,
\begin{align*}
    W_2^2(\rho,\delta_x)
    =2\bigl(1-\langle\cM_\rho,x\rangle\bigr).
\end{align*}
Thus $1-\|\cM_\rho\|_2\leq\frac12W_2^2(\rho,\delta_x)$, with equality at $x=\cM_\rho/\|\cM_\rho\|_2$ when
$\cM_\rho\neq0$. Consequently, we have
\begin{align*}
    \sqrt{2a_\nu}
    \leq W_2(\nu_T,\delta_{\cU_T^\mu})
    \leq W_2(\mu_T,\nu_T)
       +W_2(\mu_T,\delta_{\cU_T^\mu})
    =\sqrt Q+\sqrt{2a_\mu}.
\end{align*}
Squaring both sides, using $(a+b)^2 \leq 2a^2 + 2b^2$ and $Q + a_\mu \leq 1/8$, we get
\begin{equation}\label{eq: a_nu bound}
    a_{\nu} \leq Q + 2 a_\mu  \leq 1/4 .
\end{equation}
Since $R_T^\mu,R_T^\nu\geq3/4$,
\Cref{lem:finite-variation-center} gives
\begin{align*}
    \sup_{t \geq T} \| M_T^\mu - M_t^\mu \|_2 \leq (1 + \pi) a_\mu 
    , \qquad \sup_{t \geq T} \| M_T^\nu - M_t^\nu \|_2 \leq (1 + \pi) a_\nu
\end{align*}
Therefore, using \eqref{eq: a_nu bound}  and the fact that $Q + a_\mu \leq 1/8$,
\begin{equation} \label{eq: M_t bound}
\begin{split}
    \|M_t^\mu - M_t^\nu \|_2 & \leq \| M_T^\mu - M_T^\nu \|_2 + \| M_t^\mu - M_T^\mu \|_2 + \| M_t^\nu - M_T^\nu \|_2 
    \\ &\leq \sqrt{Q} + 2 (1 + \pi) (a_\mu + a_\nu) 
    \\ & \leq \sqrt{Q} + 2 (1 + \pi) (a_\mu + Q + 2 a_\mu)  
    \\ & \leq 2(2 + \pi) \sqrt{Q} + 6(1+\pi) a_\mu 
\end{split}
\end{equation}
Using \eqref{eq: a_nu bound} and \eqref{eq: M_t bound}, \eqref{eq: W2 bound} becomes 
\begin{align*}
    \sup_{t \geq T} W_2^2(\mu_t, \nu_t)  \leq C (Q + a_\mu)
\end{align*}
for $C = \max(2 + 8(2+\pi)^2, 6 + 72 (1+\pi)^2)$.
\end{proof}

\section{Propagation of chaos}
\label{sec:propagation-of-chaos}
Let $z_1(0),\ldots,z_n(0)\overset{\mathrm{i.i.d.}}{\sim}\mu_0$, and let $z_i(t)$ evolve according to \eqref{eq:particle-system-eqn}. Define the empirical measure
\begin{align*}
    \mu_t^n:=\frac1n\sum_{i=1}^n\delta_{z_i(t)}.
\end{align*}
Then $\mu_t^n$ is a weak solution of \eqref{eq:main kuramoto}
with initial datum $\mu_0^n$, so the late-time comparison in
\Cref{cor:late-time-w2} applies. To combine this comparison with
the synchronization estimate from \Cref{sec:synchronization},
we establish a finite-time propagation of chaos bound.

\subsection{Finite-time comparison}

The following one-sided estimate takes account of the relation between
the mean vector and the particles generating it.

\begin{lemma}\label{lem:one-sided-kuramoto}
For any $\rho,\nu\in\cP(\S)$ and $\pi\in\Pi(\rho,\nu)$,
\begin{align}\label{eq:one-sided-kuramoto}
    \int_{\S\times\S}
    \langle x-y,\cX_\rho(x)-\cX_\nu(y)\rangle\,\diff\pi(x,y)
    \leq\int_{\S\times\S}\|x-y\|_2^2\,\diff\pi(x,y).
\end{align}
\end{lemma}

\begin{proof}
For $x,y\in\S$ and $v,w\in\R^d$, expanding the projections gives
\begin{align}\label{eq:one-sided-projection-identity}
    \langle x-y,\proj_xv-\proj_yw\rangle
    =\langle x-y,v-w\rangle
     -\frac{\langle x,v\rangle+\langle y,w\rangle}{2}
      \|x-y\|_2^2.
\end{align}
Set
\begin{align*}
    h(x,y):=\frac{\langle x,\cM_\rho\rangle
                    +\langle y,\cM_\nu\rangle}{2},
    \qquad
    r^2:=\frac{\|\cM_\rho\|_2^2+\|\cM_\nu\|_2^2}{2}.
\end{align*}
By Cauchy--Schwarz, $|h(x,y)|\leq r$, while the marginal conditions
give $\int h\,\diff\pi=r^2$. If $r<1$, it follows that
\begin{align*}
    \int_{\S\times\S}\frac{1}{1+h}\,\diff\pi
    &\leq
    \int_{\S\times\S}\frac{1-h}{1-r^2}\,\diff\pi
    =1.
\end{align*}
If $r=1$, both $\rho$ and $\nu$ are Dirac masses, so $h=1$
$\pi$-almost surely and the same bound holds. Hence weighted
Cauchy--Schwarz yields
\begin{align*}
    \|\cM_\rho-\cM_\nu\|_2^2
    &=\left\|\int_{\S\times\S}(x-y)\,\diff\pi\right\|_2^2\\
    &\leq
      \left(\int_{\S\times\S}(1+h)\|x-y\|_2^2\,\diff\pi\right)
      \left(\int_{\S\times\S}\frac{1}{1+h}\,\diff\pi\right)\\
    &\leq\int_{\S\times\S}(1+h)\|x-y\|_2^2\,\diff\pi.
\end{align*}
Finally, integrating \eqref{eq:one-sided-projection-identity} with
$v=\cM_\rho$ and $w=\cM_\nu$ gives
\begin{align*}
    \int_{\S\times\S}
    \langle x-y,\cX_\rho(x)-\cX_\nu(y)\rangle\,\diff\pi
    =\|\cM_\rho-\cM_\nu\|_2^2
     -\int_{\S\times\S}h(x,y)\|x-y\|_2^2\,\diff\pi,
\end{align*}
which proves \eqref{eq:one-sided-kuramoto}.
\end{proof}

In particular, transporting an optimal coupling along two Kuramoto flows
gives $W_2(\rho_t,\nu_t)\leq e^{t-s}W_2(\rho_s,\nu_s)$ for $t\geq s$.

\begin{proposition}\label{prop:finite-time-w2}
There exists $C_d>0$ such that, for every $T\geq0$ and every
$\eta\in(0,1)$, with probability at least $1-\eta$,
\begin{align}\label{eq:finite-time-w2}
    \sup_{0\leq t\leq T}W_2^2(\mu_t^n,\mu_t)
    \leq C_d\left[
        \mathfrak r_{d-1}(n)
        +\sqrt{
            \frac{
                \log(2+T)+\log(2/\eta)
            }{n}}
        +\bigl(1+\log(2/\eta)\bigr)
         \frac{e^{2T}}{n}
    \right].
\end{align}
\end{proposition}

\begin{proof}
To separate the interaction error from empirical sampling, couple each particle to the mean-field characteristic starting at the same point. Using the flow from \eqref{eq:C1 diffeomorphism}, define
\begin{align}\label{eq:auxiliary-system}
    \bar z_i(t):=\phi_{0\to t}(z_i(0)),
    \qquad
    \partial_t\bar z_i(t)=\proj_{\bar z_i(t)}\cM_t,
    \qquad
    \bar z_i(0)=z_i(0).
\end{align}
Since $\mu_t=(\phi_{0\to t})_\#\mu_0$, the $\bar z_i(t)$ are i.i.d. with law $\mu_t$ for each fixed $t$. Set
\begin{align}\label{eq:auxiliary-empirical-objects}
    \bar\mu_t^n:=\frac1n\sum_{i=1}^n\delta_{\bar z_i(t)},
    \qquad
    \xi_t:=\frac1n\sum_{i=1}^n\bar z_i(t)-\cM_t,
\end{align}
and define the squared particle discrepancy
\begin{align}\label{eq:def-D}
    D_t:=\frac1n\sum_{i=1}^n\|z_i(t)-\bar z_i(t)\|_2^2.
\end{align}
The labeled pairing $z_i(t)\leftrightarrow\bar z_i(t)$ defines a coupling of $\mu_t^n$ and $\bar\mu_t^n$. Therefore
\begin{equation}\label{eq:finite-time-error-split}
\begin{aligned}
    W_2^2(\mu_t^n,\mu_t)
    &\leq 2W_2^2(\mu_t^n,\bar\mu_t^n)+2W_2^2(\bar\mu_t^n,\mu_t) \\
    &\leq 2D_t+2W_2^2(\bar\mu_t^n,\mu_t).
\end{aligned}
\end{equation}
We first bound $D_t$. The auxiliary particles solve
\begin{align*}
    \partial_t\bar z_i(t)
    =\proj_{\bar z_i(t)}
       \left(\frac1n\sum_{j=1}^n\bar z_j(t)-\xi_t\right).
\end{align*}
Apply \Cref{lem:one-sided-kuramoto} to the labeled coupling of
$\mu_t^n$ and $\bar\mu_t^n$. Cauchy--Schwarz gives
\begin{align*}
    \frac12\partial_tD_t
    =\frac1n\sum_{i=1}^n
      \left\langle z_i-\bar z_i,
      \cX_{\mu_t^n}(z_i)-\cX_{\bar\mu_t^n}(\bar z_i)
      +\proj_{\bar z_i}\xi_t\right\rangle \leq D_t+\sqrt{D_t}\,\|\xi_t\|_2.
\end{align*}
Consequently, in the sense of upper Dini derivatives,
\begin{align}\label{eq:sqrt-D-differential}
    \partial_t^+\sqrt{D_t}\leq\sqrt{D_t}+\|\xi_t\|_2.
\end{align}
Since $D_0=0$, Gronwall's inequality yields
\begin{align}\label{eq:D-convolution}
    \sqrt{D_t}
    &\leq\int_0^t e^{t-s}\|\xi_s\|_2\,\mathrm ds
    \leq e^t\mathcal Z_n,
    \qquad
    \mathcal Z_n:=\int_0^\infty e^{-s}\|\xi_s\|_2\,\mathrm ds.
\end{align}
For each fixed $t$, independence of the auxiliary particles implies
\begin{align*}
    \mathbb E\|\xi_t\|_2^2
    &=\frac1n\int_\S\|x-\cM_t\|_2^2\,\diff\mu_t(x)
     =\frac{1-R_t^2}{n}\leq\frac1n.
\end{align*}
By Cauchy--Schwarz and Fubini's theorem,
$\mathbb E\mathcal Z_n\leq1/\sqrt n$.
Moreover, replacing one initial sample changes $\xi_t$ by at most
$2/n$ at every time, and hence changes $\mathcal Z_n$ by at most
$2/n$. McDiarmid's inequality (\Cref{lem:mcdiarmid}) therefore gives
\begin{align*}
    \mathbb P\left(\mathcal Z_n>\frac1{\sqrt n}+r\right)
    \leq e^{-nr^2/2},\qquad r>0.
\end{align*}
Taking $r=\sqrt{2\log(2/\eta)/n}$ and using
\eqref{eq:D-convolution}, we obtain, with probability at least
$1-\eta/2$,
\begin{align}\label{eq:finite-time-coupling-D}
    \sup_{0\leq t\leq T}D_t
    \leq C\bigl(1+\log(2/\eta)\bigr)\frac{e^{2T}}{n}.
\end{align}
The sampling estimate in \Cref{lem:uniform-auxiliary-sampling}, also with failure probability $\eta/2$, bounds the second term in \eqref{eq:finite-time-error-split} uniformly over $[0,T]$. Combining the two bounds by a union bound proves \eqref{eq:finite-time-w2}.
\end{proof}

\subsection{Proof of the main results}\label{subsec:main-theorem-proof}

Recall the definitions in \eqref{eq:def-beta-p-epsilon}. Applying \Cref{prop:sharp exponential rate R_t} with
$\varepsilon/2$ and absorbing its linear prefactor into the remaining
exponential decay, we obtain, after enlarging the constant to cover
the initial time interval,
\begin{align}\label{eq:mean-field-rate-recalled}
    1-R_t\leq C_\varepsilon e^{-\beta_{p,\varepsilon}t},
    \qquad t\geq0.
\end{align}
The choice of $T$ below balances $e^{2T}/n$ in
\Cref{prop:finite-time-w2} with $e^{-\beta_{p,\varepsilon}T}$ in
\eqref{eq:mean-field-rate-recalled}.

\begin{proof}[Proof of \Cref{thm:uniform-time-poc-final}]
Write
\begin{align*}
    \beta:=\beta_{p,\varepsilon},
    \qquad
    \alpha:=\frac{\beta}{\beta+2},
    \qquad
    L:=1+\log(8/\eta).
\end{align*}
Suppose first that $n>L$, and set
\begin{align*}
    T:=\frac{\log(n/L)}{\beta+2}.
\end{align*}
Then
\begin{align*}
    e^{-\beta T}
    =\frac{Le^{2T}}{n}
    =\left(\frac{L}{n}\right)^\alpha.
\end{align*}
By \Cref{prop:finite-time-w2}, with probability at least $1-\eta$,
\begin{align*}
    \sup_{0\leq t\leq T}W_2^2(\mu_t^n,\mu_t)
    \leq C_d\left[
        \mathfrak r_{d-1}(n)
        +\sqrt{\frac{\log(2+T)+\log(2/\eta)}{n}}
        +\bigl(1+\log(2/\eta)\bigr)\frac{e^{2T}}{n}
    \right].
\end{align*}
The last term is at most $(L/n)^\alpha$. Since
$T\leq(\log n)/2$ and $\alpha<1/2$,
\begin{align*}
    \sqrt{\frac{\log(2+T)+\log(2/\eta)}{n}}
    \leq C_\alpha\left(\frac{L}{n}\right)^\alpha.
\end{align*}
Indeed, $\sqrt{L/n}\leq(L/n)^\alpha$, and
$\log(2+T)\leq C_\alpha n^{1-2\alpha}$ gives
\begin{align*}
    \sqrt{\frac{\log(2+T)}{n}}
    \leq C_\alpha n^{-\alpha}
    \leq C_\alpha\left(\frac{L}{n}\right)^\alpha.
\end{align*}
Also, \eqref{eq:mean-field-rate-recalled} gives
\begin{align*}
    1-R_T\leq C_\varepsilon e^{-\beta T}
    =C_\varepsilon\left(\frac{L}{n}\right)^\alpha.
\end{align*}
Consequently, \eqref{eq:proof-strategy-reduction} implies
\begin{align*}
    \sup_{t\geq0}W_2^2(\mu_t^n,\mu_t)
    \leq C\left[
        \mathfrak r_{d-1}(n)
        +\left(\frac{L}{n}\right)^\alpha
    \right].
\end{align*}
If $n\leq L$, the same bound follows from $\sup_{t\geq0}W_2^2(\mu_t^n,\mu_t)\leq4$,
since $(L/n)^\alpha\geq1$. This proves
\eqref{eq:uniform-time-poc-final}.
\end{proof}

\begin{proof}[Proof of \Cref{cor:empirical-rate-high-dimension}]
Condition
\eqref{eq:high-dimensional-condition} gives
$\beta_p/(\beta_p+2)>2/(d-1)$, so we can choose $\varepsilon>0$
with $\alpha_{p,\varepsilon}>2/(d-1)$. Also $(d-1)>4$, since $\beta_p\leq2$.
Thus $\mathfrak r_{d-1}(n)=n^{-2/(d-1)}$. If $1+\log(8/\eta) \leq n$, then
\begin{align*}
    \mathfrak r_{d-1}(n)+\left(\frac{1+\log(8/\eta)}{n}\right)^{\alpha_{p,\varepsilon}}
    \leq2\left(\frac{1+\log(8/\eta)}{n}\right)^{2/(d-1)},
\end{align*}
and the result follows from \Cref{thm:uniform-time-poc-final}.
If $1+\log(8/\eta)>n$, it follows from $\sup_{t\geq0}W_2^2(\mu_t^n,\mu_t)\leq4$.
\end{proof}

\appendix
\section{Sampling and concentration estimates}
\label{app:sampling-estimates}

\begin{lemma}[McDiarmid's inequality {\cite{mcdiarmid1989method}}]\label{lem:mcdiarmid}
Let $X_1,\ldots,X_n$ be independent, and let $F=F(X_1,\ldots,X_n)$ satisfy
\begin{align*}
    \sup_{x_1,\ldots,x_n,x_i'}
    |F(x_1,\ldots,x_i,\ldots,x_n)
      -F(x_1,\ldots,x_i',\ldots,x_n)|
    \leq c_i.
\end{align*}
Then, for every $r>0$,
\begin{align}\label{eq:mcdiarmid}
    \mathbb P(F-\mathbb EF\geq r)
    \leq\exp\left(-\frac{2r^2}{\sum_{i=1}^nc_i^2}\right).
\end{align}
The same estimate holds for $\mathbb P(\mathbb EF-F\geq r)$.
\end{lemma}

We now estimate the empirical sampling error for the auxiliary particles defined in the proof of \Cref{prop:finite-time-w2}.

\begin{lemma}\label{lem:uniform-auxiliary-sampling}
There exists $C_d>0$ such that, for every $T\geq0$ and every
$\eta\in(0,1)$, with probability at least $1-\eta$,
\begin{align}\label{eq:uniform-auxiliary-sampling}
    \sup_{0\leq t\leq T}W_2^2(\bar\mu_t^n,\mu_t)
    \leq C_d\left[
        \mathfrak r_{d-1}(n)
        +\sqrt{
          \frac{
            \log(2+T)+\log(1/\eta)
          }{n}}
    \right].
\end{align}
\end{lemma}

\begin{proof}
For every $\rho\in\cP(\S)$, the empirical transport estimate of
\cite[Theorem~1]{fournier2015rate}, applied in finitely many
bi-Lipschitz coordinate charts of the $(d-1)$-dimensional sphere, gives
\begin{align}\label{eq:fixed-time-empirical-expectation}
    \mathbb E W_2^2\left(
      \frac1n\sum_{i=1}^n\delta_{X_i},\rho
    \right)
    \leq C_d\mathfrak r_{d-1}(n),
\end{align}
where $X_1,\ldots,X_n$ are independent with common law $\rho$, and the
constant is uniform in $\rho$.

We first record the stability of the mean-field characteristic flow over
a unit time interval. Fix $s\geq0$, and let $X_t:=\phi_{s\to t}(x)$, $Y_t:=\phi_{s\to t}(y)$. By \eqref{eq:one-sided-projection-identity},
\begin{align*}
    \frac{\mathrm d}{\mathrm dt}\|X_t-Y_t\|_2 ^2
    =-\langle X_t+Y_t,\cM_t\rangle \|X_t-Y_t\|_2 ^2\leq 2\|X_t-Y_t\|_2 ^2,
\end{align*}
because $\|\cM_t \|_2 \leq1$. Hence,
\begin{align}\label{eq:auxiliary-flow-w2-stability}
    W_2^2\bigl(
      (\phi_{s\to t})_\#\rho,
      (\phi_{s\to t})_\#\nu
    \bigr)
    \leq e^{2(t-s)}W_2^2(\rho,\nu)
\end{align}
for every $t\geq s$ and every $\rho,\nu\in\cP(\S)$.

Set
\begin{align*}
    N_T:=\max\{1,\lceil T\rceil\},
    \qquad
    I_k:=[k,\min\{k+1,T\}],
    \qquad
    k=0,\ldots,N_T-1,
\end{align*}
and define
\begin{align*}
    F_k:=\sup_{t\in I_k}W_2^2(\bar\mu_t^n,\mu_t).
\end{align*}
By \eqref{eq:auxiliary-flow-w2-stability} and the identities $\bar\mu_t^n=(\phi_{k\to t})_\#\bar\mu_k^n$ and $\mu_t=(\phi_{k\to t})_\#\mu_k$, we have $F_k\leq e^2W_2^2(\bar\mu_k^n,\mu_k)$. Since the $\bar z_i(k)$ are independent with common law $\mu_k$, \eqref{eq:fixed-time-empirical-expectation} gives
\begin{align}\label{eq:block-empirical-expectation}
    \mathbb EF_k
    \leq C_d\mathfrak r_{d-1}(n),
\end{align}
after changing $C_d$.

We next apply bounded differences directly to $F_k$. Replacing one initial
sample changes only the corresponding auxiliary trajectory. Thus, at every
time $t$, the empirical measure $\bar\mu_t^n$ changes by replacing one atom.
Retaining an optimal transport plan and redirecting the mass $1/n$ emitted
by that atom to the same target submeasure shows that
\begin{align*}
    \left|
      W_2^2(\bar\mu_t^n,\mu_t)
      -W_2^2(\bar\mu_t^{n,\prime},\mu_t)
    \right|
    \leq\frac4n.
\end{align*}
Consequently, replacing one initial sample changes $F_k$ by at most $4/n$. \Cref{lem:mcdiarmid} and
\eqref{eq:block-empirical-expectation} therefore give, for every $r>0$,
\begin{align*}
    \mathbb P\left(
      F_k>C_d\mathfrak r_{d-1}(n)+r
    \right)
    \leq\exp\left(-\frac{nr^2}{8}\right).
\end{align*}
Taking a union bound over the $N_T\leq1+T$ intervals gives
\begin{align*}
    \mathbb P\left(
      \sup_{0\leq t\leq T}W_2^2(\bar\mu_t^n,\mu_t)
      >C_d\mathfrak r_{d-1}(n)+r
    \right)
    \leq N_T\exp\left(-\frac{nr^2}{8}\right).
\end{align*}
Choosing
\begin{align*}
    r=
    \sqrt{\frac{
      8\bigl(\log N_T+\log(1/\eta)\bigr)
    }{n}}
\end{align*}
and using $N_T\leq1+T$ proves
\eqref{eq:uniform-auxiliary-sampling}, after changing $C_d$.
\end{proof}

\vspace{0.1in}

\noindent {\bf Acknowledgments.} 
The authors are deeply grateful to Professor Philippe Rigollet for his generous support and insightful discussions throughout the development of this work.


\bibliographystyle{alpha}
\bibliography{references}

@incollection{mcdiarmid1989method,
  author    = {Colin McDiarmid},
  title     = {On the method of bounded differences},
  booktitle = {Surveys in Combinatorics, 1989},
  editor    = {J. Siemons},
  series    = {London Mathematical Society Lecture Note Series},
  volume    = {141},
  pages     = {148--188},
  publisher = {Cambridge University Press},
  year      = {1989},
  doi       = {10.1017/CBO9781107359949.008}
}

@misc{lin2026sharp,
  title={Sharp Convergence of {Wasserstein} Gradient Flows for Spectrally Nonnegative Interaction Energies},
  author={Lin, Zhengjiang and Rigollet, Philippe},
  note={arXiv:2609.27008},
  year={2026}
}

@article{liu2023towards,
  title={Towards understanding the dynamics of gaussian-stein variational gradient descent},
  author={Liu, Tianle and Ghosal, Promit and Balasubramanian, Krishnakumar and Pillai, Natesh},
  journal={Advances in Neural Information Processing Systems},
  volume={36},
  pages={61234--61291},
  year={2023}
}

@article{chandra2019continuous,
  title={Continuous versus discontinuous transitions in the D-dimensional generalized Kuramoto model: Odd D is different},
  author={Chandra, Sarthak and Girvan, Michelle and Ott, Edward},
  journal={Physical Review X},
  volume={9},
  number={1},
  pages={011002},
  year={2019},
  publisher={APS}
}

@article{wiesenfeld1996synchronization,
  author  = {Wiesenfeld, Kurt and Colet, Pere and Strogatz, Steven H.},
  title   = {Synchronization transitions in a disordered
             {Josephson} series array},
  journal = {Physical Review Letters},
  volume  = {76},
  number  = {3},
  pages   = {404--407},
  year    = {1996},
  doi     = {10.1103/PhysRevLett.76.404}
}

@article{kozyreff2000global,
  author  = {Kozyreff, G. and Vladimirov, A. G. and Mandel, Paul},
  title   = {Global coupling with time delay in an array
             of semiconductor lasers},
  journal = {Physical Review Letters},
  volume  = {85},
  number  = {18},
  pages   = {3809--3812},
  year    = {2000},
  doi     = {10.1103/PhysRevLett.85.3809}
}

@article{goldstein2009noise,
  author  = {Goldstein, Raymond E. and Polin, Marco and Tuval, Idan},
  title   = {Noise and synchronization in pairs
             of beating eukaryotic flagella},
  journal = {Physical Review Letters},
  volume  = {103},
  number  = {16},
  pages   = {168103},
  year    = {2009},
  doi     = {10.1103/PhysRevLett.103.168103}
}

@article{myung2015gaba,
  author  = {Myung, Jihwan and Hong, Sungho and DeWoskin, Daniel
             and {De Schutter}, Erik and Forger, Daniel B.
             and Takumi, Toru},
  title   = {{GABA}-mediated repulsive coupling between circadian
             clock neurons in the {SCN} encodes seasonal time},
  journal = {Proceedings of the National Academy of Sciences},
  volume  = {112},
  number  = {29},
  pages   = {E3920--E3929},
  year    = {2015},
  doi     = {10.1073/pnas.1421200112}
}

@inproceedings{bruno2025multiscale,
  title={A multiscale analysis of mean-field transformers in the moderate interaction regime},
  author={Bruno, Giuseppe and Pasqualotto, Federico and Agazzi, Andrea},
  booktitle={NeurIPS},
  year={2025}
}

@misc{chen2026propagation,
  title={Propagation of Chaos in Contextual Flow Maps},
  author={Shi Chen and Zhengjiang Lin and Kaizhao Liu and Philippe Rigollet},
  year={2026},
  note={arXiv:2605.16747}
}

@misc{chen2025quantitative,
      title={Quantitative Clustering in Mean-Field Transformer Models},
      author={Shi Chen and Zhengjiang Lin and Yury Polyanskiy and Philippe Rigollet},
      year={2025},
      note={arXiv:2504.14697}
}

@misc{andrew25,
      title={Synchronization of mean-field models on the circle}, 
      author={Yury Polyanskiy and Philippe Rigollet and Andrew Yao},
      year={2025},
      note={arXiv:2507.22857}
}

@misc{criscitiello2024synchronization,
      title={Synchronization on circles and spheres with nonlinear interactions},
      author={Christopher Criscitiello and Quentin Rebjock and Andrew D McRae and Nicolas Boumal},
      year={2024},
      note={arXiv:2405.18273}
}

@misc{burger2025analysis,
      title={Analysis of mean-field models arising from self-attention dynamics in transformer architectures with layer normalization},
      author={Martin Burger and Samira Kabri and Yury Korolev and Tim Roith and Lukas Weigand},
      year={2025},
      note={arXiv:2501.03096}
}

@article{geshkovski2023mathematical,
  title={A mathematical perspective on transformers},
  author={Geshkovski, Borjan and Letrouit, Cyril and Polyanskiy, Yury and Rigollet, Philippe},
  journal={Bull. Amer. Math. Soc.},
  volume={62},
  number={3},
  pages={427--479},
  year={2025},
  doi={10.1090/bull/1863}
}

@article{morales2022trend,
  title={On the trend to global equilibrium for {Kuramoto} oscillators},
  author={Morales, Javier and Poyato, David},
  journal={Annales de l'Institut Henri Poincar{\'e} C, Analyse non lin{\'e}aire},
  volume={40},
  number={3},
  pages={631--716},
  year={2023},
  doi={10.4171/AIHPC/47},
  eprint={1908.07657},
  archivePrefix={arXiv},
  primaryClass={math.AP},
  url={https://doi.org/10.4171/AIHPC/47}
}

@article{acebron2005kuramoto,
  title={The {K}uramoto model: {A} simple paradigm for synchronization phenomena},
  author={Acebr{\'o}n, Juan A and Bonilla, Luis L and P{\'e}rez Vicente, Conrad J and Ritort, F{\'e}lix and Spigler, Renato},
  journal={Reviews of modern physics},
  volume={77},
  number={1},
  pages={137--185},
  year={2005},
  publisher={APS}
}

@inproceedings{kuramoto1975self,
  title={Self-entrainment of a population of coupled non-linear oscillators},
  author={Kuramoto, Yoshiki},
  booktitle={International Symposium on Mathematical Problems in Theoretical Physics: January 23--29, 1975, Kyoto University, Kyoto/Japan},
  pages={420--422},
  year={1975},
  organization={Springer}
}

@article{geshkovski2023emergence,
  title={The emergence of clusters in self-attention dynamics},
  author={Geshkovski, Borjan and Letrouit, Cyril and Polyanskiy, Yury and Rigollet, Philippe},
  journal={Advances in Neural Information Processing Systems},
  volume={36},
  pages={57026--57037},
  year={2023}
}

@article{fournier2015rate,
  title={On the rate of convergence in Wasserstein distance of the empirical measure},
  author={Fournier, Nicolas and Guillin, Arnaud},
  journal={Probability Theory and Related Fields},
  volume={162},
  number={3-4},
  pages={707--738},
  year={2015},
  publisher={Springer}
}

@article{carrillo2011global,
 author = {J. A. Carrillo and M. DiFrancesco and A. Figalli and T. Laurent and D. Slepčev},
title = {{Global-in-time weak measure solutions and finite-time aggregation for nonlocal interaction equations}},
volume = {156},
journal = {Duke Mathematical Journal},
number = {2},
publisher = {Duke University Press},
pages = {229 -- 271},
year = {2011}
}

@inproceedings{frouvelle2019long,
  title={Long-time dynamics for a simple aggregation equation on the sphere},
  author={Frouvelle, Amic and Liu, Jian-Guo},
  booktitle={Stochastic Dynamics Out of Equilibrium: Institut Henri Poincar{\'e}, Paris, France, 2017},
  pages={457--479},
  year={2019},
  organization={Springer}
}

@inproceedings{bruno2024emergence,
title={Emergence of meta-stable clustering in mean-field transformer models},
author={Bruno, Giuseppe and Pasqualotto, Federico and Agazzi, Andrea},
booktitle={International Conference on Learning Representations},
year={2025}
}

@article{benedetto2015complete,
  title={On the complete phase synchronization for the {K}uramoto model in the mean-field limit},
  author={Benedetto, Dario and Caglioti, Emanuele and Montemagno, Umberto},
  journal={Communications in Mathematical Sciences},
  volume={13},
  number={7},
  pages={1775--1786},
  year={2015},
  publisher={International Press of Boston}
}

@article{carrillo2014contractivity,
  title={Contractivity of transport distances for the kinetic {K}uramoto equation},
  author={Carrillo, Jos{\'e} A and Choi, Young-Pil and Ha, Seung-Yeal and Kang, Moon-Jin and Kim, Yongduck},
  journal={Journal of Statistical Physics},
  volume={156},
  number={2},
  pages={395--415},
  year={2014},
  publisher={Springer}
}

@article{dobrushin1979vlasov,
  title={Vlasov equations},
  author={Dobrushin, Roland L'vovich},
  journal={Functional Analysis and Its Applications},
  volume={13},
  number={2},
  pages={115--123},
  year={1979},
  doi={10.1007/BF01077243}
}

@article{weed2019sharp,
  title={Sharp asymptotic and finite-sample rates of convergence of empirical measures in {W}asserstein distance},
  author={Weed, Jonathan and Bach, Francis},
  year={2019}, 
 volume = {25},
journal = {Bernoulli},
number = {4A},
pages = {2620--2648}
}

@article{lacker2023sharp,
  title={Sharp uniform-in-time propagation of chaos},
  author={Lacker, Daniel and Le Flem, Luc},
  journal={Probability Theory and Related Fields},
  volume={187},
  number={1-2},
  pages={443--480},
  year={2023},
  publisher={Springer},
  doi={10.1007/s00440-023-01192-x}
}

@article{rosenzweig2023global,
  title={{Global-in-time mean-field convergence for singular Riesz-type diffusive flows}},
  author={Rosenzweig, Matthew and Serfaty, Sylvia},
  journal={The Annals of Applied Probability},
  volume={33},
  number={2},
  pages={954--998},
  year={2023},
  publisher={Institute of Mathematical Statistics},
  doi={10.1214/22-AAP1833}
}

@article{gerber2025uniform,
  title={Uniform-in-time propagation of chaos for the {C}ucker--{S}male model},
  author={Gerber, Nicolai Jurek and Vaes, Urbain},
  journal={Kinetic and Related Models},
  volume={20},
  pages={56--79},
  year={2026},
  doi={10.3934/krm.2025021}
}

@article{delarue2025uniform,
  author={Delarue, Fran{\c{c}}ois and Tse, Alvin},
  title={Uniform in time weak propagation of chaos on the torus},
  journal={Annales de l'Institut Henri Poincar{\'e}, Probabilit{\'e}s et Statistiques},
  volume={61},
  number={2},
  pages={1021--1074},
  year={2025},
  doi={10.1214/23-AIHP1451}
}

@incollection{sznitman1991topics,
  author={Sznitman, Alain-Sol},
  title={Topics in propagation of chaos},
  booktitle={{\'E}cole d'{\'E}t{\'e} de Probabilit{\'e}s de Saint-Flour XIX---1989},
  editor={Hennequin, Paul-Louis},
  series={Lecture Notes in Mathematics},
  volume={1464},
  pages={165--251},
  publisher={Springer},
  address={Berlin, Heidelberg},
  year={1991},
  doi={10.1007/BFb0085169}
}

@article{malrieu2001logarithmic,
  author={Malrieu, Florent},
  title={Logarithmic {Sobolev} inequalities for some nonlinear {PDE}'s},
  journal={Stochastic Processes and their Applications},
  volume={95},
  number={1},
  pages={109--132},
  year={2001},
  doi={10.1016/S0304-4149(01)00095-3}
}

@article{delgadino2023phase,
  author={Delgadino, Mat{\'i}as G. and Gvalani, Rishabh S. and Pavliotis, Grigorios A. and Smith, Scott A.},
  title={Phase transitions, logarithmic {Sobolev} inequalities, and uniform-in-time propagation of chaos for weakly interacting diffusions},
  journal={Communications in Mathematical Physics},
  volume={401},
  number={1},
  pages={275--323},
  year={2023},
  doi={10.1007/s00220-023-04659-z},
  eprint={2112.06304},
  archivePrefix={arXiv}
}

@article{ha2018uniform,
  author={Ha, Seung-Yeal and Kim, Jeongho and Zhang, Xiongtao},
  title={Uniform stability of the {Cucker--Smale} model and its application to the mean-field limit},
  journal={Kinetic and Related Models},
  volume={11},
  number={5},
  pages={1157--1181},
  year={2018},
  doi={10.3934/krm.2018045}
}

@misc{agazzi2026quantitative,
  title={Quantitative Diffusive Limits for Singular Nonlocal Transport},
  author={Agazzi, Andrea and Bruno, Giuseppe and Pasqualotto, Federico and Rigollet, Philippe},
  year={2026},
  note={arXiv:2609.11837},
  eprint={2609.11837},
  archivePrefix={arXiv},
  primaryClass={math.AP},
  url={https://arxiv.org/abs/2609.11837}
}

@misc{glasgow2026uniform,
  author = {Glasgow, Margalit and Bruna, Joan},
  title = {Uniform-in-Time Weak Propagation-of-Chaos in Shallow Neural Networks},
  year = {2026},
  note = {arXiv:2605.22010},
  eprint = {2605.22010},
  archivePrefix = {arXiv},
  primaryClass = {stat.ML},
  url = {https://arxiv.org/abs/2605.22010}
}

@misc{gerber2025consensus,
  author = {Gerber, Nicolai and Hoffmann, Franca and Kim, Dohyeon and Vaes, Urbain},
  title = {Uniform-in-time propagation of chaos for {Consensus-Based Optimization}},
  year = {2025},
  note = {arXiv:2505.08669, revised February 2026},
  eprint = {2505.08669},
  archivePrefix = {arXiv},
  primaryClass = {math.PR},
  url = {https://arxiv.org/abs/2505.08669}
}

\end{document}